\documentclass[12pt,a4paper]{article}

\usepackage{amsmath, amssymb, verbatim, mathrsfs} 
\usepackage{theorem}

\newtheorem{theorem}{Theorem}[section] 
\newtheorem{lemma}[theorem]{Lemma} 
\newtheorem{corollary}[theorem]{Corollary}
\newtheorem{proposition}[theorem]{Proposition}
\theorembodyfont{\normalfont} 
\newtheorem{definition}[theorem]{Definition}
\newtheorem{remark}[theorem]{Remark}
  
\newcommand{\re}{\mathbb{R}} 
\newcommand{\proj}{\mathbb{P}} 
\newcommand{\cpx}{\mathbb{C}} 
\newcommand{\om}{\omega} 
\newcommand{\eps}{\varepsilon}
\newcommand{\ph}{\varphi}

\newcommand{\ra}{\rightarrow} 
\newcommand{\olo}{\mathcal{O}}

\newcommand{\vol}{\textrm{vol}}
 
\newcommand{\Ric}{\textrm{Ric}} 
 
\newcommand{\wed}{\wedge} 
\newcommand{\barr}{\overline} 
\newcommand{\del}{\partial} 
\newcommand{\delbar}{\barr{\del}}
\newcommand{\ddb}{i\del\delbar} 

\newcommand{\mab}{\mathcal{M}_{\alpha}}
\newcommand{\Bl}{\textrm{Bl}}
\numberwithin{equation}{section}   
 
\newcommand{\Aut}{\operatorname{Aut}} 

\newcommand{\dimo}[1][]         {\noindent\textbf{Proof#1}. } 

\newcommand{\fine}    {\begin{flushright} 
             \textsc{Q.E.D.} 
             \end{flushright}}

\begin{document} 
\title{Twisted cscK metrics and K\"ahler slope stability}

\author{Jacopo Stoppa} 
\date{} 
 
\maketitle 
 
\begin{abstract} 
  \noindent 
We introduce a cohomological obstruction to solving the constant scalar curvature K\"ahler (cscK) equation twisted by a semipositive form, appearing in works of Fine and Song-Tian.

Geometrically this gives an obstruction for a manifold to be the base of a holomorphic submersion carrying a cscK metric in certain ``adiabatic'' classes. In turn this produces many new examples of general type threefolds with classes which do not admit a cscK representative. 

When the twist vanishes our obstruction extends the slope stability of Ross-Thomas to effective divisors on a K\"ahler manifold. Thus we find examples of non-projective slope unstable manifolds.
\end{abstract} 
\section{Introduction}
One of the central open problems in K\"ahler geometry is to characterise the K\"ahler classes represented by a constant scalar curvature K\"ahler (cscK) form. 

In the the algebraic case of representing the first Chern class of a positive line bundle this is related to stability in algebraic geometry by a conjecture of Yau \cite{yau}, Tian \cite{ding_tian} and Donaldson \cite{don_toric}. 

It is a fundamental result in the theory that the existence of a cscK metric in the algebraic case implies K-semistability in the sense of Donaldson (the shortest proof is given in \cite{don_calabi}). 

From this Ross-Thomas derived a cohomological obstruction to solving the cscK equation in the algebraic case known as slope stability for polarised manifolds \cite{ross_thomas}. This proved to be a very effective tool in the study of cscK metrics on projective bundles and algebraic surfaces, see e.g. \cite{ross_panov}.

In this paper we expand from this starting point in two directions.\\
 
Firstly with Theorem \ref{richard_conj} we extend slope stability (with respect to effective divisors) to any K\"ahler class on a K\"ahler manifold, confirming a conjecture of Ross-Thomas. This is explained in \ref{ross_thomas} below. Non-projective examples are given in \ref{non_proj_example}.

Our methods are necessarily different from those of Ross-Thomas, in particular they are differential-geometric in nature. The key fact in this connection is the lower bound for the K-energy proved by Chen-Tian, see \ref{k_energy} below. Then Section \ref{stability_proof} is essentially devoted to showing that the slope stability condition with respect to an effective divisor arises precisely by imposing that this lower bound holds as the K\"ahler form concentrates along the divisor to some extent dictated by positivity.\\  

The other main theme of this paper is that K\"ahler slope stability generalises to an equation introduced by Fine and Song-Tian and which we call the twisted cscK equation. Solutions to the twisted cscK equation arise as zeros of the moment map for an action of the group of exact symplectomorphisms. This was pointed out by G. Sz\'ekelyhidi in analogy to the case of cscK metrics where it was shown by Donaldson \cite{don_fields} and Fujiki. The proof is given in Section \ref{moment_map}.

Thus our main obstruction result, Theorem \ref{stability_theorem}, is stated in this more general setup.

This extenstion is not a mere formality however, as obstructing the twisted cscK equation leads to interesting geometric applications through the so-called ``adiabatic limit" construction.

We hope to make this clear in the rest of this Introduction, see in particular the obstruction Theorem \ref{adiabatic_theorem} and the application to general type threefolds explained in \ref{threefolds} below.

\subsection{The twisted cscK equation.}
Let $M$ be a compact K\"ahler manifold with K\"ahler class $\Omega$ and let $\alpha$ be a closed pointwise semi-positive $(1,1)$-form on $M$. The equation that we are interested in, and which we call the \emph{twisted cscK equation} is finding a metric $\omega\in\Omega$ such that
\begin{equation}\label{twist}
 S(\omega) - \Lambda_\omega\alpha = \widehat{S}_\alpha.
\end{equation}
Here $S(\omega)$ is the scalar curvature, $\Lambda_{\om}\alpha$ denotes the trace of $\alpha$ with respect to $\om$ 
\[\Lambda_{\om}\alpha\,\om^n = n\alpha \wed \om^{n-1}\]
and $\widehat{S}_\alpha$ is the only possible topological constant, given by
\[ \widehat{S}_\alpha =
\frac{n(c_1(X)-[\alpha])\cup[\omega]^{n-1}}{[\omega]^n}. \]
In particular if $\alpha=0$ we recover the cscK equation,
\begin{equation}\label{csc}
S(\om) = \widehat{S},
\end{equation}
and $\widehat{S}$ is just the average scalar curvature. 

Equation \ref{twist} is a generalisation of \ref{csc} which arises naturally in the work of Fine \cite{fine_surfaces}, \cite{fine_fibrations} and Song-Tian \cite{song_tian}. In these cases twisted cscK metrics in dimension $n$ are related to genuine cscK metrics in higher dimension by a limiting process.

Let $\pi: M \ra B$ be a holomorphic submersion of K\"ahler manifolds endowed with a relatively ample line bundle $L$. The so-called \textit{adiabatic classes} on $M$ are given by
\begin{equation}\label{adiabatic}
\Omega_r = c_1(L) + r\pi^* \Omega_B
\end{equation}
where $\Omega_B$ is any K\"ahler class on $B$. These are certainly K\"ahler for large $r$. 

Suppose that \ref{csc} can be solved in all the fibres of $\pi$ inside the restriction of $c_1(L)$, and that the fibres and base carry no nontrivial holomorphic vector fields. In this case Fine proved that \ref{csc} can be solved in $\Omega_r$ for $r \gg 0$ provided that the twisted equation \ref{twist} is solvable in $\Omega_B$ with respect to a special choice of $\alpha$. Loosely speaking the relevant $\alpha$ is the pullback of the Weil-Petersson type metric from the moduli space of cscK metrics on a fibre. 

We take up the application to adiabatic cscK metrics in Section \ref{adiabatic_classes}. In particular Theorem \ref{adiabatic_theorem} gives a necessary numerical condition for their existence.

Recently Song-Tian studied the \emph{untraced} form of \ref{twist} (see also \cite{fine_fibrations}), 
\begin{equation}\label{genKE}
\Ric(\om) = \lambda\om + \alpha
\end{equation}
in connection with the K\"ahler-Ricci flow on an elliptic surface. In this case $\lambda = -1$ and $\alpha$ is given by the pullback of the genuine Weil-Petersson metric plus a singular contribution corresponding to multiple fibres. 

They prove that the K\"ahler-Ricci flow converges to a solution of \ref{genKE} on the base, and call this a \emph{generalised K\"ahler-Einstein metric}. 

Similarly Song-Tian introduce the equation
\[\delbar\,\nabla^{(1,0)}(S(\om) - \Lambda_{\om}\alpha) = 0\]
generalising the \textit{extremal equation} of Calabi. This specialises to \ref{twist} when there are no holomorphic vector fields, and to \ref{genKE} when 
\[c_1(M)- [\alpha] = \lambda \Omega.\]
We refer to loc. cit. Section 5 for a complete discussion.  
\subsection{Obstructions.} Our first result gives a cohomological obstruction for the twisted cscK equation to admit a solution in the class $\Omega$. 

Recall that the \emph{positive cone} 
\[\mathcal{P}\subset H^{1,1}(M, \cpx) \cap H^2(M, \re)\]
is given by cohomology classes $\Theta$ which evaluate positively on irreducible subvarieties, that is 
\[\int_V \Theta^p > 0\]
for any $p$-dimensional irreducible analytic subvariety $V \subset M$. 

Let $D \subset M$ be an effective divisor.
\begin{definition}
The \emph{Seshadri constant} of $D$ with respect to $\Omega$ is given by
\[\epsilon(D, \Omega) = \sup \{x : \Omega - x c_1(\olo(D)) \in \mathcal{P}\}.\]
\end{definition}
We define coefficients $\alpha_i$, $i = 0, 1$ by
\begin{align}
&\alpha_1(x) = \frac{(\Omega-x c_1(D))^n}{n!},\\
&\alpha_2(x) = \frac{(c_1(M)-[\alpha])\cup(\Omega-x c_1(D))^{n-1}}{2(n-1)!}.
\end{align}
\begin{definition} The \emph{twisted Ross-Thomas polynomial} of $(X, \Omega)$ with respect to $D$ and $\alpha$ is given by 
\[\mathcal{F}_{\alpha}(\lambda) = \int^{\lambda}_{0}(\lambda-x)\,\alpha_2(x)dx + \frac{\lambda}{2}\alpha_1(0)-\frac{\widehat{S}_{\alpha}}{2}\int^{\lambda}_0 (\lambda-x)\,\alpha_1(x)dx.\]
\end{definition} 
\begin{theorem}[Stability condition]\label{stability_theorem}
If \ref{twist} is solvable in the class $\Omega$ then 
\begin{equation}\label{stability_condition}
\mathcal{F}_{\alpha}(\lambda) \geq 0
\end{equation}
for all effective divisors $D \subset X$ and $0 \leq \lambda \leq \epsilon(D, \Omega)$.
\end{theorem}
The proof will be given in Section \ref{stability_proof}.
\begin{remark}
For $\lambda < \eps(D,\Omega)$ we have $\alpha_1(x) > 0$ so the quotient
\[\mu_{\lambda}(\olo_D, \Omega) = \frac{\int^{\lambda}_{0}(\lambda-x)\,\alpha_2(x)dx + \frac{\lambda}{2}\alpha_1(0)}{\int^{\lambda}_0 (\lambda-x)\,\alpha_1(x)dx}\]
is well defined. Following Ross-Thomas we call this \emph{slope}. Then the inequality \ref{stability_condition} can be rewritten as
\begin{equation}\label{stability_conditionII}
\mu_{\lambda}(\olo_D, \Omega) \geq \frac{\widehat{S}_{\alpha}}{2}.
\end{equation}
The point is that choosing $\Omega = c_1(L)$, $\alpha = 0$ gives back the slope stability condition of Ross-Thomas \cite{ross_thomas} in the algebraic setting. 

The relation of our work with slope stability is explained in \ref{ross_thomas} below.
\end{remark}

Next we obtain an obstruction to the existence of the adiabatic cscK holomorphic submersions introduced above.
\begin{theorem}\label{adiabatic_theorem}
Let $\pi: M \ra B$ be a holomorphic submersion endowed with a relatively ample line bundle $L$, which is fibrewise cscK with fibrewise average scalar curvature $S_b$. If the adiabatic classes \ref{adiabatic} admit cscK metrics for $r \gg 0$ then \ref{stability_condition} holds with
\begin{equation}
[\alpha] = c_1(\pi_* K_{X|B}) + \frac{S_b}{n+1}c_1(\pi_*L)
\end{equation}
for all effective divisors $D \subset B$ and $0 \leq \lambda \leq \epsilon(D, \Omega_B)$.
\end{theorem}
An analogous statement holds when $L$ is replaced by a class $\Omega_0$ which is positive and cscK along the fibres. The proof will be given in Section \ref{adiabatic_classes}.\\

We will give concrete examples of how one can apply both results in Section \ref{examples}. In particular Corollary \ref{surfaces} gives an a priori obstruction for an algebraic surface with a fixed K\"ahler class to be the base of an adiabatic cscK submersion.
\subsection{K\"ahler slope stability.}\label{ross_thomas} Ross-Thomas prove that subschemes can give an obstruction to the solvability of the cscK equation in the projective case $\Omega = c_1(L)$. Given a subscheme $Z \subset M$ they define a cohomological slope $\mu_{\lambda}(\olo_Z, L)$ depending on a positive parameter $\lambda$ less than the Seshadri constant. They prove that if \ref{csc} is solvable in $c_1(L)$ then the slope inequality 
\begin{equation}\label{slope_inequality}
\mu_{\lambda}(\olo_Z,L) \geq \mu(M, L)
\end{equation} 
where $\mu(M, L) = \frac{-n K_M.L^{n-1}}{2L^n}$ must hold for all $Z \subset M$ and $0 < \lambda \leq \epsilon(Z, L)$. In other words, $(M,L)$ must be \textit{slope semistable}. The particular form of \ref{slope_inequality} is of course modelled on slope stability for vector bundles (in terms of quotient sheaves).  

As we already observed in the special case of effective divisors \ref{slope_inequality} is equivalent to inequality \ref{stability_condition} with the choices $\Omega = c_1(L)$, $\alpha = 0$. 

Motivated by the corresponding situation for vector bundles Ross-Thomas extended the definition of slope to any analytic subvariety of a K\"ahler manifold, and conjectured that \ref{slope_inequality} still gives a necessary condition for \ref{csc} to be solvable. We refer to \cite{ross_thomas} Section 4.4 for more details. In our case we only need to allow any K\"ahler class $\Omega$ in place of $c_1(L)$. 

Thus Theorem \ref{stability_theorem} with the choice $\alpha = 0$ gives a proof of this conjecture in the case of effective divisors, namely  
\begin{theorem}\label{richard_conj}
If a K\"ahler class $\Omega$ on a K\"ahler manifold $M$ admits a cscK representative it is slope semistable with respect to effective divisors. 
\end{theorem}  
A concrete non-projective example is given in Lemma \ref{non_proj_example}.
\begin{remark}
The case of divisors is somehow central in slope stability. In general blowing up the subscheme $Z \subset M$ reduces to the exceptional divisor $E \subset \Bl_Z M$, polarised by $L^{1/\lambda}\otimes\olo(-E)$. This trick is not directly relevant to us because for our methods we need smoothness of the ambient space. However one might hope to keep track of the resolution of singularities in this context. Since all our concrete examples are based on effective divisors, we do not pursue this here. 

In a different vein Ross-Panov \cite{ross_panov} Section 4 prove that a slope unstable \emph{surface} is destabilised by a divisor. 
\end{remark}
\subsection{General type 3-folds.}\label{threefolds} By the Aubin-Calabi-Yau Theorem if $K_M$ is ample we can find a K\"ahler-Einstein metric in its first Chern class. On the other hand one may ask what happens for classes far from $c_1(K_M)$. These have proved hard to obstruct and the first example of a surface with ample canonical bundle and a class with no cscK representative was found only recently by Ross \cite{ross_inv}. 

To date the only examples for general type 3-folds are somewhat trivial, namely products of an obstructed surface with a curve of genus at least 2, endowed with the product class. In \ref{threefold_exampleI}, \ref{threefold_exampleII}, \ref{threefold_exampleIII} we apply our results on adiabatic cscK metrics to provide new  examples of general type 3-folds with obstructed classes.
\subsection{Log-geometry.} Given an effective divisor $D$ on $M$, log-geometry replaces the canonical bundle $K_M$ by $K_M + D$ and calls $D$ the \emph{boundary} of $M$. Note that this is precisely what happens when we replace slope stability with our condition \ref{stability_conditionII}, that is we replace $c_1(K_M)$ by $c_1(K_M) + [\alpha]$ throughout. On the other hand Sz\'ekelyhidi \cite{gabor} has conjectured that the slope stability condition for the existence of a complete extremal metric on the complement of a reduced effective divisor $D$ is precisely log-slope, replacing $c_1(K_M)$ by $c_1(K_M + \olo(D))$ in the slope inequality. Thus in geometric terms coupling with an adiabatic class in a holomorphic submersion and removing a divisor should have the same slope stability condition.
\subsection{K-energy.}\label{k_energy} Finally a few words about our method of proof. Our results are based on the computation of a suitable energy functional, the ``twisted K-energy", along a metric degeneration breaking the K\"ahler form into currents of integration (see Section \ref{stability_proof} and Theorem \ref{limit} in particular). 

Thus we introduce the natural analogue of the K-energy on the space of K\"ahler metrics in a given K\"ahler class. For any $\phi \in \mathcal{H}$, the space of K\"ahler potentials of $\om$, let $\om_{\phi} = \om + \ddb\phi$, $\mu_{\phi} = \frac{\om_{\phi}^n}{n!}$.
\begin{definition}The \emph{variation of the twisted K-energy} at $\om_{\phi}$ is the 1-form  
\[ \delta\mathcal{M}_\alpha(\delta \phi) =
-\int_M\delta\phi(S(\omega_{\phi})-\Lambda_\omega\alpha-\widehat{S}_\alpha)\,\mu_{\phi}\]
for $\delta \phi \in T_{\phi}\mathcal{H} \cong C^{\infty}(M, \re)$.
\end{definition}

The form $\delta\mab$ on $\mathcal{H}$ is closed. Choosing a base-point $\om \in \Omega$ and integrating along any path gives a well-defined \emph{twisted K-energy} $\mab$.

This is just the sum of the well known K-energy \cite{mabuchi} and a variant of the J-functional introduced by Chen \cite{chen_mabuchi}. Its critical points are precisely solutions of \ref{twist}. It is not hard to prove that these are actually local minimisers.  

Since it is the sum of two functionals which are convex along geodesics in $\mathcal{H}$, $\mathcal{M}_\alpha$ is also convex. More precisely the second derivative along a path $\phi = \phi_t$ is given by
\begin{equation}
\frac{d^2}{dt^2}\mathcal{M}_{\alpha} = \|\delbar \nabla^{(1,0)}\dot{\phi}\|^2_{\phi} + (\del \dot{\phi} \wed \delbar \dot{\phi}, \alpha)_{\phi} - \int_{M}(\ddot{\phi} - \frac{1}{2}|\nabla^{1,0}\dot{\phi}|^2_{\phi} )(S(\om_{\phi})-\Lambda_{\om_{\phi}}\alpha - \widehat{S}_{\alpha})\,\mu_{\phi}
\end{equation}
where all the metric quantities are computed with respect to $\om + \ddb \phi$. The \emph{geodesic equation} is $\ddot{\phi} - \frac{1}{2}|\nabla^{1,0}\dot{\phi}|^2_{\phi} = 0$. Semmes and Donaldson have shown that this is the equation of geodesics for a negatively curved locally symmetric space structure on $\mathcal{H}$. The existence of geodesic segments is well known to be equivalent to a boundary value problem for a degenerate complex Monge-Amp\'ere equation.

If either $\alpha$ is strictly positive at a point or $M$ has no Hamiltonian holomorphic vector fields, $\mab$ is \textit{strictly} convex along geodesics in $\mathcal{H}$. In this case, the deep results of Chen-Tian \cite{chen_tianII} on the regularity of weak geodesics imply uniqueness of solutions of \ref{twist}. In general, they imply the following lower bound on $\mab$.
\begin{lemma}\label{lo_bound}
If \ref{twist} is solvable the twisted K-energy $\mab$ is bounded below in $\Omega$.
\end{lemma}
\textbf{Acknowledgements}. This work took shape mainly thanks to conversations with G. Sz\'ekelyhidi. The moment map interpretation in particular is due to him. I thank my advisor R. Thomas for his great support. He also first suggested to prove Theorem \ref{richard_conj} via K-energy. J. Fine, A. Ghigi, D. Panov, S. Rollenske and J. Ross provided useful comments and suggestions. The support of the Geometry Groups at Imperial College and Pavia, in particular G. Pirola, has been essential. Finally I am grateful to S. Donaldson for many interesting discussions about cscK metrics and the material presented here. 
\section{Moment map interpretation}\label{moment_map}
In this section we show that the operator
\[\om \mapsto \widehat{S}_{\alpha}-S(\om) + \Lambda_{\om}\alpha\]
can be viewed as a moment map, essentially by combining Donaldson's moment map computations in \cite{don_fields} and \cite{don_moment}. 

Let $(B,\omega)$ be a symplectic manifold of dimension $n$ and let us
assume for simplicity that $H^1(B)=0$. We write $\mathscr{J}$ for the space of integrable 
complex structures on $B$ which are compactible with $\omega$. The space $\mathscr{J}$ admits a natural symplectic form, and the action of the group of symplectomorphisms $\mathrm{Symp}(B,\omega)$ preserves this symplectic form. Let us identify the Lie algebra of $\mathrm{Symp}(B,\omega)$ with $C^\infty_0(B,\re)$, the smooth functions with vanishing integral, via the Hamiltonian construction. Also, using
the $L^2$ inner product with respect to the volume form $\frac{\omega^n}{n!}$ let
us identify the dual of $C^\infty_0(B,\re)$ with itself. 
\begin{theorem}[Donaldson] The map
	\[ J \mapsto \hat{S}-S(J),\]
  is an equivariant moment map for the action of
  $\mathrm{Symp}(B,\omega)$ on $\mathscr{J}$. Here $S(J)$ is the scalar
  curvature of the K\"ahler metric given by $(\omega, J)$ and $\hat{S}$
  is its average which is independent of $J$. 
\end{theorem}

Now let $M$ be diffeomorphic to $B$ and fix a symplectic form $\alpha$
on $M$. Let $\mathscr{M}$
be the space of diffeomorphisms $f:B\to M$ in a fixed homotopy class. 
The tangent space
$T_f\mathscr{M}$ to $\mathscr{M}$ at $f$ can be identified with the
space of vector fields on $M$. This has a natural symplectic form given
by 
\[ \Omega_f(v,w) = \int_B f^*(\alpha(v,w)) \frac{\omega^n}{n!}.\]
The group $\mathrm{Symp}(B,\omega)$ acts on $\mathscr{M}$ by composition
on the right, preserving $\Omega$. 
\begin{theorem}[Donaldson] The map
  \[ \nu : f \mapsto \Lambda_\omega f^*(\alpha)\]
  is an equivariant moment map for the action of
  $\mathrm{Symp}(B,\omega)$ on $\mathscr{M}$. 
\end{theorem}
\dimo We give the simple proof for the reader's convenience.
Let us first consider the infinitesimal action of
$\mathrm{Symp}(B,\omega)$ on $\mathscr{M}$. For this let $H\in
C^\infty(B)$, and write $X_H$ for the Hamiltonian vector field on $B$
generated by $H$. The infinitesimal action at a point $f\in\mathscr{M}$ is
given by $f_*(X_H)$. To show that $\nu$ is a moment map we therefore
have to show
\[ \langle d\nu_f(w), H\rangle = \Omega_f(f_*X_H, w) =
\int_S f^*(\alpha(f_*X_H, w)) \frac{\omega^n}{n!}, \]
where $w\in T_f\mathscr{M}$. 

Let $f_t$ be a path in $\mathscr{M}$ such that $f_0=f$ and
$\frac{d}{dt}f_t=w$ at $t=0$. Then
\[ 
  d\nu_f(w) = \left.\frac{d}{dt}\right|_{t=0} \Lambda_\omega
  f^*(\alpha) = \Lambda_\omega f^*(d(\iota_w\alpha)).
\]
Therefore 
\[ \begin{aligned}
  \langle d\nu_f(w), H\rangle &= \int_B H d(f^*(\iota_w\alpha)) \wedge
  \frac{\omega^{n-1}}{(n-1)!} \\
  &= \int_B f^*(\iota_w\alpha) \wedge\iota_{X_H}\omega \wedge
  \frac{\omega^{n-1}}{(n-1)!}\\
  &= \int_B f^*(\alpha(f_*X_H, w)) \frac{\omega^n}{n!}. 
\end{aligned} \]
This shows that $\nu$ is a moment map.\fine

We now combine the above two moment maps. 
Let us fix a complex structure $I$ on $M$ which is compatible with
$\alpha$ and consider the space
\[ \mathscr{S} = \{(f,f^*(I))\,|\, f\in\mathscr{M}\text{ such that
}f^*(I)\text{ is compatible with }\omega\}\subset\mathscr{M}\times
\mathscr{J}.\]
Then $\mathscr{S}$ is preserved by the action of
$\mathrm{Symp}(B,\omega)$, and the moment map restricted to
$\mathscr{S}$ is given by 
\[ (f,f^*(I)) \mapsto -S(\omega,f^*(I)) + 
\Lambda_\omega f^*(\alpha) + C, \]
where $C$ is a constant such that the integral over $B$ is zero. 

Following Donaldson let us now consider the complexification of the
action of $\mathrm{Symp}(B,\omega)$. This can only be done on the level
of the Lie algebras. Both $\mathscr{J}$ and $\mathscr{M}$ have natural
complex structures so we can complexify the infinitesimal action of the
Lie algebra $C^\infty_0(B)$. This gives rise to an integrable
distribution on $\mathscr{M}\times\mathscr{J}$ which is tangent to
$\mathscr{S}$. We think of the leaves of the resulting foliation of
$\mathscr{S}$ as the complexified orbits. 

If $(f,f^*(I))$ is a zero of the moment map, then 
\[ S( (f^{-1})^*\omega, I) - \Lambda_{(f^{-1})^*\omega} \alpha \] 
is constant. If in addition $(f,f^*(I))$ is 
in the complexified orbit of $(\mathrm{id},
I)$ then in fact $(f^{-1})^*\omega$ is in the same K\"ahler class as
$\omega$, so we have a solution of the twisted cscK equation. 
\section{Proof of stability condition}\label{stability_proof}
As already mentioned in the Introduction what we actually prove in this section is a result on the asymptotic behaviour of the twisted K-energy.
\begin{theorem}\label{limit}
Let $D \subset M$ be an effective divisor. For any base-point $\om \in \Omega$ and $0 < \lambda < \epsilon(D, \Omega)$ there exists a family of K\"ahler forms $\om_{\eps} \in \Omega,\,0 < \eps \leq 1$ with $\om_1 = \om$ and such that as $\eps \ra 0$ 
\begin{equation}\label{limit_equ}
\mab(\eps) = -\pi\mathcal{F}_{\alpha}(\lambda)\log(\eps) + O(1)
\end{equation} 
\end{theorem}
By Lemma \ref{lo_bound} this settles Theorem \ref{stability_theorem}. 
\begin{remark}
It must be pointed out that expansions for the K-energy of the form \ref{limit_equ} are known to hold for the Fubini-Study metric along the fibres of a $\cpx^*$-action in projective space under some regularity assumptions, thanks to the work of Paul-Tian \cite{paul_tian}, Phong-Sturm-Ross \cite{p_s_ross} and others. 
\end{remark}
Theorem \ref{limit} will be proved in several steps. We first recall the fundamental Nakai-Moishezon criterion of Demailly-Paun \cite{dem_pau}.
\begin{theorem}[Demailly-Paun] The K\"ahler cone $\mathcal{K}$ is a connected component of the positive cone $\mathcal{P}$.
\end{theorem} 
\begin{corollary}\label{dem_pau}
\[\epsilon(D, \Omega) = \sup \{x: \Omega-x c_1(\olo(D)) \in \mathcal{K}\}.\]
\end{corollary}
This holds because $\Omega-x c_1(\olo(D))$ is in the connected component of $\Omega \in \mathcal{K}$. 

Let $\om \in \Omega$ be any K\"ahler form, and $h$ be any Hermitian metric on $\olo(D)$ with curvature form $\Theta$. We begin by defining the family in Theorem \ref{limit} near $\eps = 0$. 

Given a canonical section $\sigma \in H^0(\olo(D))$ we define potentials 
\[\psi_{\eps} = \frac{1}{2}\log(\eps^2 + |\sigma|^2_h)\]
for $0 < \eps \ll 1$.
By the Poincar\'e-Lelong equation there is a weak convergence 
\[i\del\delbar\lambda\psi_{\eps} \rightharpoonup -\lambda \Theta + \lambda[D]\]
where $[D]$ denotes the (closed, positive) current of integration along $D$. By Corollary \ref{dem_pau} for $0 < \lambda < \epsilon(D, \Omega)$ we can find a potential $u$ (independent of $\eps$) such that 
\[\eta = \om -\lambda \Theta + \ddb u > 0.\]
We define our family for $0 < \eps \ll 1$ by
\begin{equation}\label{sequence}
\om_{\eps} = \om + \ddb u + \ddb \lambda\psi_{\eps}.
\end{equation}
Since $[D]$ is positive,  
\[\eta + \lambda [D] > \eta\]
holds in the sense of currents, and the sequence $\om_{\eps}$ converges weakly to  $\eta + \lambda[D]$, so the inequality
\[\om_{\eps} > \eta \]
holds in the sense of currents for $\eps \ll 1$. In other words $\om_{\eps} - \eta $ is a strictly positive $(1,1)$-current for all small $\eps$. But since it is also a smooth $(1,1)$-form, this implies it is actually a K\"ahler form (more generally a similar pointwise statement holds for a current with locally $L^1$ coefficients, see \cite{dem_gaz} Section 3). 

We conclude that for $0 < \eps \ll 1$, $\om_{\eps}$ is a K\"ahler form. We still need to prescribe the base-point to be $\om_1 = \om$. This can be achieved by choosing $\psi_{\eps} \ra 0$ as $\eps \ra 1$, and making $u = u_{\eps}$ dependent on $\eps$ away from $\eps = 0$ so that $u_{\eps} \ra 0$ as $\eps \ra 1$ and $\om + \ddb \psi_{\eps} + \ddb u_{\eps}$ is always K\"ahler. 

The next important observation is that we have a uniform $C^{\infty}$ bound for $\om_{\eps}$ away from $D$ (that is, in any compact of $X \setminus D$). This means that for the sake of proving \ref{limit} we need only study the behaviour of $\mab$ near points of supp$(D)$. 

Let us first write down in detail the model case of $x \in \textrm{supp}(D)$ near which $D$ is reduced and smooth. This will be enough for the general case thanks to standard results in the theory of currents.

So choose coordinates $(z = \{z_i\}^{n-1}_{i = 1}, w)$ near $x$ such that $w$ is a local generator for $\olo(D)$. Then
\[\psi_{\eps} = \frac{1}{2}\log(\eps^2 + e^{-2\ph}|w|^2)\]  
near $x$, where $e^{-\ph}$ is the weight for $h$ on $\olo(D)|_U$ with respect to the Euclidean norm $|\cdot|$. We make the first order expansion in the $D$-transversal direction
\[e^{-2\ph} = e^{-2\widetilde{\ph}}(1+O(|w|))\]
where $\widetilde{\ph} = \widetilde{\ph}(z)$ does not depend on the transversal coordinate $w$. Moreover from now on we assume that we have chosen normal $z$-coordinates at $x$ so that 
\[\widetilde{\ph} = \del_{z_i}\widetilde{\ph} = \del_{\barr{z}_j}\widetilde{\ph} = 0.\]
Since $\widetilde{\ph}$ is independent of $w$ this holds in a small slice $\{z = z(x)\}$. 

We introduce a piece of notation that will be very useful for the rest of this section: we will denote by $r$ any smooth function $r = r(z, w)$ which is $O(|w|)$ uniformly in $U$, or any smooth form whose coefficient functions have the same property in our fixed coordinates.

With this convention, direct calculation gives
\begin{lemma}
In the slice $\{z = z(x)\}$,     
\begin{align*}
&\del_{\eps}\psi_{\eps} = \eps(\eps^2 + (1+r)|w|^2)^{-1};\\
&\del_{w}\del_{\barr{w}} \psi_{\eps} = \frac{1}{2}(1+r)\,\eps^2(\eps^2 + (1+r))|w|^2)^{-2};\\
&\del_{z_i}\del_{\barr{w}} \psi_{\eps} = r;\\
&\del_{z_i}\del_{\barr{z}_j} \psi_{\eps} = -(1+r)\,\del_{z_i}\del_{\barr{z}_j}\widetilde{\ph}\, |w|^2( \eps^2 + (1+r))|w|^2)^{-1}.
\end{align*}
\end{lemma}
Taking exterior powers we find
\begin{lemma}\label{powers} In the slice $\{z = z(x)\}$ for $p \geq 1$
\begin{equation*}
(\ddb \psi_{\eps})^p = r + (-1)^p(1+r)^p|w|^{2 p}(\eps^2 + (1+r)|w|^2)^{-p}(\ddb\widetilde{\ph})^p +
\end{equation*}
\begin{equation*}
\frac{(-1)^{p-1}}{2}(1+r)^{p-1}\eps^2|w|^{2(p-1)}(\eps^2 + (1+r)|w|^2)^{-p-1} p\,(\ddb \widetilde{\ph})^{p-1}\wed i\,dw\wed d\barr{w},
\end{equation*}
and so
\begin{equation*}
\del_{\eps}\psi_{\eps}(\ddb \psi_{\eps})^p = r + (-1)^p(1+r)^p\eps |w|^{2p}(\eps^2 + (1+r)|w|^2)^{-p-1}(\ddb \widetilde{\ph})^p + 
\end{equation*}
\begin{equation*}
\frac{(-1)^{p-1}}{2}(1+r)^{p-1}\eps^3|w|^{2(p-1)}(\eps^2 + (1+r)|w|^{2})^{-p-2} p\,(\ddb\widetilde{\ph})^{p-1}\wed i\,dw\wed d\barr{w}.
\end{equation*}
\end{lemma}
This first order expansion is used to prove a global result about the weak convergence of some quantities which we will need later on when computing the K-energy.
\begin{proposition}\label{currents} The following hold in the sense of currents:
\begin{equation}\label{curr_1}
\eps\del_{\eps}\psi_{\eps}\rightharpoonup 0
\end{equation}
and for $p \geq 1$
\begin{equation}\label{curr_2}
\eps\del_{\eps}\psi_{\eps}(\ddb \psi_{\eps})^p \rightharpoonup \frac{(-1)^{p+1}\pi}{2 (1+p)}\Theta^{p-1}\wed[D].
\end{equation}
\end{proposition}
\dimo Weak convergence can be checked locally. Consider first the model case of a point $x \in \textrm{supp}(D)$ near which $D$ is reduced and smooth. Choosing coordinates at $x$ as in Lemma \ref{powers} we find that
\[\eps \del_{\eps}\psi_{\eps} = \eps^2(\eps^2 + (1+r)|w|^2)^{-1}.\]
The right hand side is uniformly bounded and converges to 0 uniformly away from $w = 0$, thus it converges to 0 weakly as $\eps \ra 0$.\\

The case when $D$ is smooth but not reduced near $x$ can be handled by a finite ramified cover $w \mapsto w^m$, where $m$ is the local multiplicity. The above weak convegence is unaffected. Thus \ref{curr_1} holds away from the set of singular points of supp$(D)$.

Weak convergence extends to all $M$ thanks to the Support and Skoda theorems, as in the proof of Lemma 2.1 of \cite{dem_pau} for example.\\

As for \ref{curr_2} multiplying by $\eps$ the second expansion in \ref{powers} we find
\begin{equation*}
\eps\del_{\eps}\psi_{\eps}(\ddb \psi_{\eps})^p = \eps\,r + (-1)^p(1+r)^p\eps^2 |w|^{2p}(\eps^2 + (1+r)|w|^2)^{-p-1}(\ddb \widetilde{\ph})^p + 
\end{equation*}
\begin{equation*}
\frac{(-1)^{p-1}}{2}(1+r)^{p-1}\eps^4|w|^{2(p-1)}(\eps^2 + (1+r)|w|^{2})^{-p-2} p\,\Theta^{p-1}\wed i\,dw\wed d\barr{w}.
\end{equation*}
The sequence of functions
\begin{equation*}
\eps\,r + (-1)^p(1+r)^p\eps^2 |w|^{2p}(\eps^2 + (1+r)|w|^2)^{-p-1} 
\end{equation*}
is uniformly bounded as $\eps \ra 0$ an converges to 0 away from $\{w = 0\}$, thus it gives no contribution to the weak limit.\\ 

The curvature form $\Theta$ appears since $\ddb \widetilde{\ph}$ represents the $w$-constant extension to $U$ of the pullback of the curvature form $\Theta$ to $D$.\\

Therefore \ref{curr_2} holds in a neighborhood of $x$ provided the sequence of forms
\[f_{\eps} = (1+r)^{p-1}\eps^4|w|^{2(p-1)}(\eps^2 + (1+r)|w|^{2})^{-p-2}dw\wed d\barr{w}\]
converges to $\frac{\pi}{p\,(p+1)}\delta_{\{w = 0\}}dw\wed d\barr{w}$ in the sense of currents. 

For this note that $f_{\eps}$ is converging uniformly to $0$ away from $w = 0$. Moreover the change of variable $w = \eps w'$ shows that the transversal integrals \[\int_{\{z = z(x)\}} f_{\eps} i\,dw\wed d\barr{w}\]
are converging to 
\[2\pi\int^{\infty}_0 \frac{s^{2p-1}}{(1+s^2)^{p+2}}ds = \frac{\pi}{p(p+1)}\]
by Lebesgue's dominated convergence.

Global weak convergence then follows as for \ref{curr_1}.
\fine
We will make use of a well known integration by parts formula due to Chen and Tian (see e.g. \cite{chen_mabuchi}), adapted to our situation.
\begin{lemma}\label{parts}
\[\delta\mab = \delta \int_M \log\left(\frac{\om^n_{\phi}}{\om^n}\right)\frac{\om^n_{\phi}}{n!} + \widehat{S}_{\alpha}\delta I + \delta J_{\alpha},\]
where
\[\delta I = \int_M \delta \phi \frac{\om^n_{\phi}}{n!},\]
\[\delta J_{\alpha} = \int_M \delta \phi\,(\emph{Ric}(\om)-\alpha)\wed\frac{\om^{n-1}_{\phi}}{(n-1)!}.\]
\end{lemma}
With the help of this formula the proof of Theorem \ref{limit} can be split into three separate statements.

We first prove an expansion for the $J_{\alpha}$-functional.
\begin{definition}
\[\mathcal{F}_{J_{\alpha}}(\lambda) = \sum^{n-1}_{p = 1} \frac{(-1)^{p+1}\lambda^{p+1}}{2(p+1)!(n-1-p)!} \int_M([\alpha]-c_1(X))\cup\Omega^{n-1-p}\cup c_1^p(\olo(D)).\]
\end{definition}  
\begin{lemma}\label{jay}
\[J_{\alpha}(\eps) = -\pi\mathcal{F}_{J_{\alpha}}(\lambda)\log(\eps) + O(1).\]
\end{lemma}
\dimo We will actually find the limit of $\eps \del_{\eps}J_{\alpha}$ as $\eps \ra 0$. The integrand for this functional is 
\[(n-1)!^{-1}\eps\del_{\eps}\lambda\psi_{\eps}(\Ric(\om)-\alpha)\wed\om^{n-1}_{\eps}\]
which by the binomial theorem applied to $(1,1)$ forms can be rewritten as
\[(n-1)!^{-1}(\Ric(\om)-\alpha)\wed \sum^{n-1}_{p = 0} \binom{n-1}{p} (\om + \ddb u)^{n-1-p}\wed \lambda^{p+1}\eps\del_{\eps}\psi_{\eps}(\ddb \psi_{\eps})^p.\]
By Corollary \ref{currents} this converges weakly to
\[\pi(n-1)!^{-1}\sum^{n-1}_{p = 1} \binom{n-1}{p} \frac{(-1)^{p+1}\lambda^{p+1}}{2(p+1)}(\Ric(\om)-\alpha)\wed (\om + \ddb u)^{n-1-p}\wed\Theta^{p-1}\wed[D].\]
Integrating over $M$ proves our claim. \fine
An identical argument applies to the $I$ functional.
\begin{definition}
\[\mathcal{F}_I(\lambda) = \sum^n_{p = 1}\frac{(-1)^{p+1}\lambda^{p+1}}{2(p+1)!(n-p)!}\int_M \Omega^{n-p}\cup c^p_1(\olo(D)).\]
\end{definition}
\begin{lemma}
\[I(\eps) = \pi\mathcal{F}_I(\lambda) \log(\eps) + O(1).\]
\end{lemma}

Finally we consider the $\int\log\det$ functional. This requires a slighlty better estimate than for the $I$ and $J_{\alpha}$ functionals. 
\begin{definition}
\[\mathcal{F}_{\log} = \sum^n_{p = 1} \frac{(-1)^{p-1}\lambda^p}{p!(n-p)!}\int_M\Omega^{n-p}\cup c_1(\olo(D))^p.\]
\end{definition}
\begin{lemma}
\[\int_M \log\left(\frac{\om^n_{\eps}}{\om^n}\right)\mu_{\eps} = -\pi\mathcal{F}_{\log} \log(\eps) + O(1).\]
\end{lemma}
\dimo We study the weak limit as $\eps \ra 0$ of the sequence of forms on $M$
\begin{equation*}
\log(\eps)^{-1}\log\left(\frac{\om^n_{\eps}}{\om^n}\right)\mu_{\eps}.
\end{equation*}
Weak convergence can be checked locally, so we can assume 
\[\log(\eps)^{-1}\log\left(\frac{\om^n_{\eps}}{\om^n}\right)\mu_{\eps} = f_{\eps} - f'_{\eps}\]
\begin{equation*} 
f_{\eps} = (\log(\eps)^{-1}\log(\det g_{\eps})\det g_{\eps})\,dz\wed d\barr{z}\wed dw\wed d\barr{w},
\end{equation*}
\begin{equation*} 
f'_{\eps} = (\log(\eps)^{-1}\log(\det g)\det g_{\eps})\,dz\wed d\barr{z}\wed dw\wed d\barr{w}.
\end{equation*}
Our first claim is that the sequence of forms $f_{\eps}$ converges weakly to 0. As in the proof of \ref{currents} we work in a good slice $\{z = z(x)\}$. The sequence certainly converges uniformly to $0$ away from $w = 0$. The claim follows if we can show that the integrals $\int f_{\eps}$ go to $0$ as $\eps \ra 0$. We will actually show that the integrals in the $w$ direction $\int_{\{z = z(x)\}}\log(\det g_{\eps})\det g_{\eps}\,dw\wed d\barr{w}$ are uniformly bounded. We can prove this using our previous computation in Lemma \ref{powers}, namely there is a decomposition
\begin{equation}\label{volume}
\om^n_{\eps} = \alpha_{\eps} + \beta_{\eps}
\end{equation}
where
\begin{equation*}
\alpha_{\eps} = r + \sum^n_{p = 0}\binom{n}{p}(-1)^p(1+r)^p\lambda^p |w|^{2p}(\eps^2+(1+r)|w|^2)^{-p}(\om + \ddb u)^{n-p}\wed (\ddb\widetilde{\ph})^p,
\end{equation*}
\begin{align*}
\beta_{\eps} = \sum^n_{p = 1}\binom{n}{p} & (\om + \ddb u)^{n-p}\wed \\ & \wed\frac{(-1)^{p-1}}{2}(1+r)^{p-1} \lambda^p \eps^2|w|^{2(p-1)} (\eps^2 + (1+r)|w|^2)^{-p-1}\times \\ & p\,(\ddb\widetilde{\ph})^{p-1}i\,dw\wed d\barr{w}.
\end{align*} 
We make the change of variable $w = \eps w'$. Under this change of variable $\alpha_{\eps}$ and $\beta_{\eps}$ pull back to forms defined transversally on $0 \leq s = |w'| < \eps$. Moreover the pullback of $\alpha_{\eps}$ is uniformly $O(\eps^2)$, while the pullback of $\beta_{\eps}$ restricted to each slice $\{z = z(x)\}$ is dominated by (in angular coordinates) 
\[C\cdot\frac{s^{2p-1}}{(1+s^2)^{p+1}}\,ds\,d\theta = O(s^{-3})\,ds\,d\theta\] 
since $p \geq 1$ (where $C$ is a positive uniform constant).

It follows that $\int \log(\det g_{\eps})\det g_{\eps}\,dw\wed d\barr{w}$ is asymptotic in each slice to 
\[\int^{\eps^{-1}}_0 \log(O(\eps^2) + O(s^{-3}))(O(\eps^2) + O(s^{-3}))\,ds.\]
By dominated convergence these integrals are uniformly bounded as $\eps \ra 0$.

Thus $f_{\eps}$ converges to $0$ weakly as $\eps \ra 0$.\\

Let us compute the weak limit of $f'_{\eps}$. By the decomposition \ref{volume}, 
\[f'_{\eps} = \log(\eps)^{-1}\log(\det(g))\,\alpha_{\eps} + \log(\eps)^{-1}\log(\det(g))\,\beta_{\eps}.\]
The sequence of forms $\log(\det(g))\,\alpha_{\eps}$ is uniformly bounded in $\eps$, therefore $\log(\eps)^{-1}\log(\det(g))\,\alpha_{\eps}$ converges weakly to $0$. 

On the other hand making the change of variable $w = \eps w'$ pulls back the sequence of forms $\log(\eps)^{-1}\log(\det(g))\,\beta_{\eps}$ to 
\[\log(\eps)^{-1}(\log(\eps^{2}) + \log(g(w')))\,\beta_{\eps}(w') = 2\,\beta_{\eps}(w') + \log(\eps)^{-1}\log(g(w'))\,\beta_{\eps}(w').\]
Now $\log(\eps)^{-1}\log(g(w'))\,\beta_{\eps}(w')$ converges weakly to $0$, while $2\,\beta_{\eps}$ converges weakly to 
\[\pi \sum^n_{p = 1} \frac{(-1)^{p-1}\lambda^p}{2 p!(n-p)!}(\om + \ddb u)^{n-p}\wed (\ddb \widetilde{\ph})^{p-1}\wed [D].\]
This is proved exactly like in Proposition \ref{currents} and Lemma \ref{jay}; in other words by the definition of $\beta_{\eps}$ it is enough to prove that for $p \geq 1$ the sequence of forms
\[(1+r)^{p-1} \lambda^p \eps^2|w|^{2(p-1)} (\eps^2 + (1+r)|w|^2)^{-p-1}dw\wed d\barr{w}\]
converges weakly to $\frac{\pi}{p}\delta_{\{w = 0\}}dw\wed d\barr{w}$. These forms are certainly converging to $0$ uniformly away from $\{w = 0\}$. Moreover the integrals in the $w$ direction of their pullbacks under $w = \eps w'$ converge to
\[2\pi\int^{\infty}_0 \frac{s^{2p-1}}{(1+s^2)^{2p+1}}ds = \frac{\pi}{p}\] 
by dominated convergence.

Integration over $M$ gives the required global contribution $-\mathcal{F}_{\log}$.
\fine
The proof of \ref{limit} is completed by the following combinatorial identity. 
\begin{lemma} 
\[\int^{\lambda}_{0}(\lambda-x)\alpha_2(x)dx + \frac{\lambda}{2}\alpha_1(0) = \mathcal{F}_{\log}(\lambda) + \mathcal{F}_{J_{\alpha}}(\lambda);\]
\[\int^{\lambda}_0 (\lambda-x)\alpha_1(x)dx  = -2\mathcal{F}_I.\]
\end{lemma}
\begin{remark}
Our proof of Theorem \ref{limit} is a metric analogue of the \textit{deformation to the normal cone} used by Ross-Thomas. In other words we compute the twisted K-energy along a sequence of metrics on $M$ which converge to the Fubini-Study metric on the fibres of $\proj(\nu(D)\oplus\olo)\ra D$
locally after rescaling (at least near smooth, reduced points of $D$). 
\end{remark}
\section{The K-energy of a holomorphic submersion}\label{adiabatic_classes}
Let $\pi: M \ra B$ be a holomorphic submersion. For the rest of this section $n$ denotes the relative dimension, while $m$ is the dimension of the base. 

We suppose there is a class $\Omega_0 \in H^{1,1}(M, \cpx)\cap H^2(M, \re)$ whose restriction to any fibre $M_b$ is a K\"ahler class containing a cscK metric $\om_b$, depending smoothly on $b$. Fitting together these forms gives a representative $\om_0 \in \Omega_0$ whose fibrewise restriction is cscK. 

The fibrewise scalar curvature and volume are fixed constants $S_b$, $\vol(M_b)$ respectively.

Choosing a K\"ahler form $\om_B \in \Omega_B$ on the base, we are interested in the existence of a cscK metric in the adiabatic classes $\Omega_r = \Omega_0 + r\pi^*\Omega_B$, $r \gg 0$. 

We describe the aforementioned result of Fine in more detail. Note that the forms $\om_b$ give a metric on the vertical tangent bundle $V$ and so on the line $\det(V)$. Taking the curvature we get a $(1,1)$-form $F_V$ representing $c_1(V)$. We denote by $F_{VH}$ the purely horizontal component of $F_V$ with respect to the horizontal-vertical decomposition induced by $\om_0$. Define a $(1,1)$-form $a$ by $F_{VH} = -\vol(X_b) a$. The choice of sign is to follow our convention that $\alpha$ (as defined below) is a \textit{semipositive} rather then \textit{seminegative} form. Taking the fibrewise average of $-F_{VH}$ with respect to $\om_0$ yields a $(1,1)$-form $\alpha$ on the base, in other words
\begin{align*}
\alpha_b = \pi_* (a\,\om^n_0) = \int_{X_b} a \,\om^n_b = \frac{\int_{X_b} -F_{VH} \om^n_b}{\vol(X_b)}
\end{align*}
where the integrals are form-valued. Semipositivity of $\alpha$ is granted by the following result. See the Introduction to \cite{fine_fibrations} for a detailed discussion.
\begin{theorem}[Fujiki-Schumacher \cite{fuj_schu}, Fine \cite{fine_fibrations}] The form $\alpha$ is semipositive. If the fibres have no nontrivial holomorphic vector fields and the submersion does not induce an isotrivial fibration on any curve in $B$ the form $\alpha$ is strictly positive.
\end{theorem}
Fine assumes that that the fibres $M_b$ have no nontrivial holomorphic vector fields and that the equation
\begin{equation}\label{joel}
S(\om_B) - \Lambda_{\om_B}\alpha = \widehat{S}_{\alpha}  
\end{equation}
is solvable for a K\"ahler metric $\om_B \in \Omega_B$ on the base. Moreover $\om_B$ must admit no nontrivial cohomologous deformations which are still solutions. The conclusion is the existence of a cscK metric in all adiabatic classes $\Omega_r$ for $r \gg 0$. 
\begin{remark}
It follows from our discussion of the K-energy in Section 1 that the condition about deformations of $\om_B$ is certainly satisfied if either $B$ has no nontrivial holomorphic vector fields or $\alpha$ is strictly positive at a single point. 
\end{remark}
Note that it is \emph{not} known if solving \ref{joel} is necessary for the existence of adiabatic cscK metrics. We expect this to be a very difficult question. Instead, we will use the K-energy to provide an obstruction.\\
  
The variation of the twisted K-energy in this case is best written in the form
\[\delta \mab (\phi) = -\int_B \delta \phi (S(\om_{B, \phi}) - \Lambda_{\om_{B,\phi}}\alpha - \widehat{S}_{\alpha})\frac{\om^n_{B,\phi}}{n!}\]
where $\om_{B,\phi} = \om_B + \ddb \phi$. 

We will show that $\vol(M_b)\mab$ is the leading order term when we expand the genuine K-energy $\mathcal{M}$ on $M$ in the adiabatic limit. Thus if the adiabatic classes admit a cscK representative the twisted K-energy on the base must be bounded below. In turn Theorem \ref{limit} gives an obstruction to the existence of these metrics. More precisely, let 
\[\om_r = \om_0 + r\om_B.\]
We are concerned with the variation of the K-energy on $M$ with respect to 
\[\delta \om_r = \delta(\om_0 + r\om_{B, \phi}) = r\,\ddb\delta \phi,\]
that is  
\[\delta\mathcal{M} = -\int_M r\delta{\phi}(S_r-\widehat{S_r})\om^{n+m}_r.\]
\begin{lemma}\label{K_fibre}
\[\delta\mathcal{M} = \emph{vol}(M_b)\delta \mab r^{m} + O(r^{m-1}).\]
\end{lemma}
\begin{corollary}\label{K_fibreII}
If the class $\Omega_r$ admits a cscK representative for $r \gg 0$ then the twisted K-energy $\mab$ is bounded below in the class $\Omega_B$ on the base.
\end{corollary}
We need a preliminary computation of scalar curvature. This is implicit in \cite{fine_surfaces} Theorem 8.1, but we need to write it down in full in order to compute the relevant K-energy.

Introduce the vertical Laplacian $\Delta_V$, characterised by
\[\Delta_V u\, \om^{n}_{b} = (\ddb u)_{VV},\]
the purely vertical component. 
\begin{lemma}
Then
\[S(\om_r) = S(\om_b) + r^{-1}(S(\om_B) - \emph{vol}(M_b)\Lambda_{\om_B} a + \Delta_V(\Lambda_{\om_B}\om_H)) + O(r^{-2}).\]
\end{lemma}
\dimo The K\"ahler form $\om_0$ on $M$ gives rise to a horizontal-vertical decomposition of forms; in particular
\[\om_0 = \om_{b} + \om_H.\]
We need to compute the Ricci form of $\om_r$. By the exact sequence of holomorphic vector bundles 
\[0\ra V \ra TM \ra H \ra 0\]
this is the sum of the curvatures of the induced metrics on the line bundles $\det(V)$, $\det(H)$, say $F_V$, $F_H$. By definition of $a$
\[F_V = \Ric(\om_b) - \vol(X_b) a\]
For $\det(H)$, we get 
\[F_H = \ddb\log(\om_H + r\om_{B})^m = \ddb\log \om^m_{B}(1 + r^{-1}\frac{\om_H\wed\om^{m-1}_B}{\om^m_B} + O(r^{-2}))\]
\[= \Ric(\om_B) + \ddb\log(1+ r^{-1}\Lambda_{\om_B}\om_H + O(r^{-2}))\]
\[= \Ric(\om_B) + r^{-1}\ddb\Lambda_{\om_B}\om_H + O(r^{-2}).\]
Taking traces,
\[\Lambda_{\om_r}F_V = S(\om_b) + \frac{\vol(F) a \wed (\om_H + r\om_B)^{m-1}}{(\om_H + r\om_B)^{m}}\]
\[= S(\om_b) - \frac{r^{m-1}\vol(M_b)a \wed \om^{m-1}_B + O(r^{m-2})}{r^m\om^m_B + O(r^{m-1})} = S(\om_b) - r^{-1}\vol(M_b)\Lambda_{\om_B}a + O(r^{-2});\]
\[\Lambda_{\om_r}F_H = \frac{\Ric(\om_B) \wed (\om_H + r\om_B)^{m-1}}{(\om_H + r\om_B)^{m}} + r^{-1}\Delta_V \Lambda_{\om_B}\om_H+\frac{\ddb(r^{-1}\Lambda_{\om_B}\om_H) \wed (\om_H + r\om_B)^{m-1}}{(\om_H + r\om_B)^{m}}\]
\[= r^{-1}(S(\om_B) + \Delta_V \Lambda_{\om_B}\om_H) + O(r^{-2}).\]
\fine
\dimo[ of Lemma \ref{K_fibre}] Note that
\[\pi_* \om^{n+m}_r = \vol(M_b) (\om_H + r\om_B)^m = \vol(M_b)r^m\om_B^m + O(r^{m-1}),\]
and
\[-\vol(M_b)\pi_* \Lambda_{\om_B} a\,\om^{n+m}_r = -\vol(M_b)\Lambda_{\om_B}\alpha\,r^m\om_B^m + O(r^{m-1}),\]
\[\pi_* \Delta_V \Lambda_{\om_B}(\om_H)\om^{n+m}_r = 0.\]
As a consequence
\[\pi_* S_r\,\om^{m+n}_r\]
\[=\left(S_b + r^{-1}\left(S(\om_B) - \Lambda_{\om_B}\alpha \right)\right)\vol(M_b)r^m\om^m_B + O(r^{m-2}).\]
The average of $S_r$ can be computed using this pushforward. We get 
\[S_b + r^{-1}\left(\widehat{S(\om_B)} - \vol(B)^{-1}\int_B \Lambda_{\om_B}\alpha\,\om^m_{B}\right) + O(r^{-2})\]
\[= S_b + r^{-1}\widehat{S_{\om_B}} + r^{-1}\vol(B)^{-1}\int_B \Lambda_{\om_B}
\alpha\,\om^m_{B} + O(r^{-2})\]
\[= S_b + r^{-1}\widehat{S}_{\alpha} + O(r^{-2}).\]
So we see that for the K-energy
\[\delta\mathcal{M} = -\int_B \pi_* (r\delta\phi (S_r - \widehat{S_r})\om^{n+m}_r)\]
\[= \int_B (r\delta \phi)\, r^{-1}(S(\om_{B, \phi}) - \Lambda_{\om_{B,\phi}}\alpha - \widehat{S}_{\alpha})(\vol(M_b)r^{m}\om^{m}_{B,t} + O(r^{m-1}))\]
\[= \vol(M_b)\delta\mab\,r^{m} + O(r^{m-1}).\]
\fine
We are finally in a position to prove our obstruction to adiabatic cscK submersions, Theorem \ref{adiabatic_theorem}.\\

\dimo[ of Theorem \ref{adiabatic_theorem}] By Corollary \ref{K_fibreII} the result follows immediately from the asymptotics of the twisted K-energy, Theorem \ref{limit}, once we know the cohomology class of $\alpha$. This is computed in the Lemma below. 
\fine
\begin{lemma}
Suppose $\Omega_0 = c_1(L)$ for a relatively ample line bundle $L$. Then
\[[\alpha] = c_1(\pi_* K_{M|B}) + \frac{S_b}{n+1}c_1(\pi_*L).\]
\end{lemma}
\dimo By definition
\[[\alpha] = [\pi_*(a\,\om^n_0)] = \vol(M_b)^{-1}[\pi_*(\Ric(\om_b)\wed\om^n_0) - \pi_*(F_V\wed\om^n_0)]\]
\[= \frac{S_b}{n+1} \pi_*c_1(L) + \pi_*(c_1(K_{M|B})),\]
and the pushforwards $\pi_* L, \pi_* K_{M|B}$ are locally free. \fine
\section{Examples}\label{examples}
\subsection{General type 3-folds.} We begin with a simple adaptation of a result of Ross-Panov \cite{ross_panov}. 

Recall that an effective divisor $D = \sum_i m_i D_i$ ($D_i$ irreducible components) is \emph{exceptional} if $[D_i.D_j]$ is negative definite. Let $p_a(D)$ be the arithmetic genus of $D$, given by $\chi(\olo_D) = 1 - p_a(D)$. 
\begin{lemma}\label{ross_panov}
Let $M$ be an algebraic surface containing an exceptional divisor $D$, $\alpha$ a semi-positive form. If the inequality  
\begin{equation}
2p_a(D)-2 + D.[\alpha] > 0
\end{equation}
holds, there exists a K\"ahler class $\Omega$ for which $\mathcal{F}_{\alpha}(\lambda) < 0$ for all small positive $\lambda$.
\end{lemma}
\dimo This is essentially the argument in loc. cit. Theorem 3.2. We adapt this to our case for the sake of completeness. 

On surfaces the stability condition \ref{stability_condition} in its slope formulation is best rewritten as
\begin{equation}\label{slope_exc1}
\frac{3(\Omega.D -\lambda((K_M + [\alpha]).D + D^2))}{2\lambda(3\Omega.D-\lambda D^2)} \geq \frac{-(K_X+[\alpha]).\Omega}{\Omega^2}
\end{equation}   
when $\Omega - \lambda D$ is a positive class. If we were allowed to choose $\Omega$ so that $\Omega.D = 0$, by adjunction we could simplify this inequality to 
\begin{equation}\label{slope_exc2}
\frac{3(2 p_a(D)-2 + D.[\alpha])}{2\lambda D^2} \geq \frac{-(K_X+[\alpha]).\Omega}{\Omega^2}.
\end{equation}   
Since $D^2 < 0$, whenever $2 p_a(D)-2 + D.[\alpha] > 0$ the left hand side diverges to $-\infty$ as $\lambda \ra 0$, thus violating \ref{stability_condition}.

Of course this argument is not quite rigorous since $\Omega$ is positive. What we do instead is to construct a sequence $\Omega_s$ of K\"ahler classes degenerating to a \emph{big} class $\Omega_0$ such that $\Omega_0.D = 0$. Then \ref{slope_exc2} must hold by continuity, provided we also have a strictly positive lower bound on the Seshadri constants $\epsilon(\Omega_s, D)$ in order to be able to substitute some positive value for $\lambda$.

Let $K$ be any reference K\"ahler form. We claim there is a choice of strictly positive numbers $r_i$ such that the sequence
\[\Omega_s = (1+s)K + \sum_i r_i D_i\]
satisfies our requests. 

The symmetric matrix $[D_i.D_j]$ is negative definite, so by diagonalising we see that we can find strictly positive numbers $r_i$ with 
\[\sum_i D_i D_j r_i = -K.D_j\]
for all $j$. This settles $\Omega_0.D = 0$. 

The lower bound on $\epsilon(\Omega_s, D)$ is more delicate. If we make the extra assumption
\[D.D_i \leq 0\]
for all $i$, a lower bound is simply given by $\min_i\{\frac{r_i}{m_i}\} > 0$. 
One can prove that the extra assumption causes no loss of generality for algebraic surfaces, see loc. cit. Corollary 3.4 for details. 
\begin{remark}
This complication only arises if $D$ has more than one irreducible component. We can avoid it in any example when $D$ can be chosen irreducible, but not in our examples later on, since a key ingredient is to take finite coverings. 
\end{remark}
Finally, $\Omega^2_0 = K.\Omega_0 \geq K^2 > 0$, so $\Omega_0$ is big. \fine
\begin{corollary}\label{surfaces}
If an algebraic surface $M$ contains an exceptional divisor $D$ with arithmetic genus $p_a(D)$ $\geq 2$ then there exists a class $\Omega$ for which $\mab$ is unbounded below for all $\alpha \geq 0$. If $p_a(D) \geq 1$ the same holds for all $\alpha > 0$. 

In the first case $(M, \Omega)$ cannot be the base of an adiabatic cscK submersion.
\end{corollary}
We will show how this can be applied to construct classes which do not admit a cscK representative on a 3-fold of general type. The idea is to start from a Kodaira surface $S$ with a suitable exceptional divisor $D$. By a result of Morita we obtain a 3-fold holomorphic submersion $X \ra S'$ upon taking a finite covering $S' \ra S$.
\begin{theorem}[Morita \cite{morita}]\label{morita} Let $S$ be a surface with a holomorphic submersion to a smooth curve, with base and fibres of genus at least 2. There exists a finite covering $S' \ra S$ and a non-isotrivial holomorphic submersion $\mathcal{X} \ra S'$ whose fibres are curves of genus at least 2.
\end{theorem}
\begin{remark}
Fine \cite{fine_fibrations}, Section 4 uses Morita's result to give new examples of cscK metrics. We will do the converse, using it to obstruct some classes.
\end{remark}
Then we use Corollary \ref{surfaces} and Corollary \ref{K_fibreII} to obstruct the adiabatic classes on $X$. We will do a few cases in detail.
  
\subsubsection{Product base} Let $C$ be a smooth curve of genus $g$ at least 2. Consider the product $S = C \times C$ with diagonal $\delta$ and fibre-classes $f_1, f_2$. $S$ has ample canonical bundle $K_S = (2 g - 2)(f_1 + f_2)$ and $\delta$ is an exceptional divisor (of genus genus $g$) since $\delta^2 = 2-2g$. 

If $S'\ra S$ is any finite covering then $K_{S'} > 0$ and the preimage $\delta'$ is again an exceptional divisor of genus at least 2. 

So choose a finite covering $S' \ra S$ as in \ref{morita} to get a 3-fold $X \ra S'$. Fix any K\"ahler class $H$ on $S'$. Define K\"ahler classes $\Omega_s = (1+s)H + \sum_i r_i \delta_i$ as in the proof of \ref{ross_panov}. 
\begin{lemma}\label{threefold_exampleI}
The general type 3-fold $X$ admits no cscK metrics in the classes 
\[a\,c_1(K_{X|S'}) + r \pi^* \Omega_s\]
for $a \in \re^{+}$, $s \ll 1$, $r \gg 0$. 
\end{lemma}
\begin{remark} The original surface example of Ross \cite{ross_inv} is in fact the self-product of a nongeneric curve $C$.
\end{remark}
\subsubsection{Atiyah-Hirzebruch base} Let $C$ be a smooth curve with genus $g$ at least $2$ and a free action of a group $G$ of order $d$ on it. This provides a finite $d$-fold covering $\pi: B \ra C$. Let $\Sigma \subset B \times C$ be the union of the graphs of $g\circ\pi$, $g \in G$. It is well known that there exists a $d$-fold branched covering $\pi_1: S \ra B \times C$, branched precisely along $\Sigma$, and that $S$ is smooth. The canonical bundle $K_S$ is ample, and the pre-image $\Gamma_1$ of the graph $\Gamma \subset B \times C$ of $\pi$ is an exceptional divisor of genus at least two, since $\Gamma^2 = d(2-2g)$. Moreover, the natural map $S \ra C$ is a holomorphic submersion. This means that we can apply Morita's construction: there is a finite covering $\pi_2: S' \ra S$ such that $S'$ is the base of a non-isotrivial holomorphic submersion $Y \ra S'$ whose fibres are smooth curves of genus at least 2. Again $K_{S'} > 0$ and the preimage $\Gamma_2$ of $\Gamma_1$ under $\pi_2$ is an exceptional divisor of genus at least 2. Setting $\Omega_s = (1+s)H + \sum_i r_i \Gamma_{2,i}$ as before (for any K\"ahler class $H$ on $S'$) we find
\begin{lemma}\label{threefold_exampleII}
The general type 3-fold $Y$ admits no cscK metrics in the classes
\[a\,c_1(K_{Y|S'}) + r \pi^* \Omega_s\]
for $a\in\mathbb{R}^+$, $s \ll 1$, $r \gg 0$. 
\end{lemma}
\begin{remark} Before the general results of \cite{ross_panov} obstructed classes on Atiyah-Hirzebruch surfaces had been constructed by Shu \cite{shu}. 
\end{remark}
\subsubsection{Catanese-Rollenske base} Another class of examples can be derived from a recent construction of Catanese-Rollenske \cite{sonke}. Let $C$ be a smooth curve of genus $g$ at least $2$ and choose a nonempty subset $\mathcal{S}\subset\Aut(C)$ (note that $\mathcal{S}$ is not necessarily a subgroup). An element $s \in \mathcal{S}$ determines the graph $\Gamma_s \subset S = C\times C$. Assume that $\Gamma_s \cap \Gamma_t = \emptyset$ for $s \neq t$; this translates into a group-theoretic condition in $\Aut(C)$ which holds for many choices. This disjoint union gives a divisor $D = \cup_{s \in \mathcal{S}}\Gamma_s \subset S$. Choose a fibre $F \subset S$ for the projection $p_1$, and fix the topological type of a ramified cover $F'\ra F$, that is a surjection $\psi: \pi_1(F\setminus D)\ra \pi_1(F) \ra 0$. Catanese-Rollenske prove that after a base change induced by a finite covering $f: C' \ra C$, there exists a ramified cover $r: S' \ra f^*S$, ramified exactly along $f^*D$, which is a Kodaira fibration for the projection $p_1\circ f\circ r$, and such that the class of the induced ramified cover of fibres is the prescribed $\psi$. 

Now we take one more finite covering $\rho: S'' \ra S'$ in order to apply Morita's Theorem and get a Kodaira fibration $Z \ra S''$ with smooth curve fibres of genus at least 2. For all $t\in \mathcal{S}$, $\Gamma''_t = \rho^* f^* r^* \Gamma_t$ is an exceptional divisor of genus at least 2. Write $\Gamma''_{t,i}$ for the irreducible components. For any K\"ahler class $H$ on $S'$ we form the classes $\Omega_s = (1+s)H + \sum_i r_i \Gamma''_{t,i}$ as in the previous examples.
\begin{lemma}\label{threefold_exampleIII}
The general type 3-fold $Z$ admits no cscK metrics in the classes
\[a\,c_1(K_{Z|S''}) + r \pi^* \Omega_s\]
for $a \in \re^{+}$, $s \ll 1$, $r \gg 0$. 
\end{lemma} 
\subsection{Slope unstable non-projective manifolds}\label{non_proj} One might try to prove Theorem \ref{richard_conj} by a deformation argument. For example it is a classical result of Kodaira that any compact K\"ahler surface admits arbitrarily small deformations which are projective. If we could perturb the cscK metric at the same time this would give an alternative proof of \ref{richard_conj} for surfaces. 

This motivates us to look for a genuinely non-projective example in higher dimensions. This is based on Voisin's manifold, combined with the following result for projective bundles. 
\begin{lemma}\label{bundles}
Let $E \ra B$ be a vector bundle on the K\"ahler manifold $(B, \Omega_B)$. Let $\olo_{\proj}(1)$ denote the relative hyperplane line bundle. If $E$ is Mumford-destabilised by a corank 1 subbundle $F \subset E$ then $\proj(E)$ admits no cscK metrics in the classes 
\[\Omega_r = c_1(\olo_{\proj}(1)) + r \pi^* \Omega_B\]
for $r \gg 0$.
\end{lemma} 
\dimo This is the K\"ahler version of a special case of \cite{ross_thomas} Theorem 5.12. 

Since $\proj(F)\subset \proj(E)$ is a divisor, all the computions in loc. cit. hold in $NE^1(\proj(E))$ and therefore carry over without any change to our situation, with the care of replacing the line bundle $L^{\otimes r}$ there with the form $r \Omega_B$. Note in passing that a large part of the proof in loc. cit. is devoted to prove that $\epsilon(\proj(F), \Omega_r) = 1$ for $r \geq r_0$ where $r_0$ does not depend on $F$ (using boundedness of quotients). This is required to establish a general correspondence with Mumford-stability. However for our weak statement we only need to prove $\epsilon(\proj(F), \Omega_r) = 1$ for $r \geq r_0(F)$. This is trivial since $c_1(\olo_E(1))-x\,c_1(\olo_F(1))$ is relatively ample for $0 < x < 1$. \fine
Voisin \cite{voisin} constructed K\"ahler manifolds in all dimensions $\geq 4$ which are not homotopy equivalent to projective ones. The simplest example in dimension $4$ is obtained by a torus $T$ with an endomorphism $f: T \ra T$. Consider the product $T \times T$ with projections $p_i$, $i = 1, 2$. There are four sub-tori in $T \times T$ given by the factors $T_i = p_i^* T$, the diagonal $T_3$ and the graph $T_4$ of $f$. Blow up the intersections of all $T_k$, that is a finite set $Q \subset T \times T$. Then blow up once more along the proper transforms of the $T_k$'s. The result of this process is of course a K\"ahler manifold $M$. Voisin proves that for a special choice of $(T, f)$, $M$ is not homotopy equivalent to (and so a fortiori not deformable to) a projective manifold. 

For any $q \in Q$, let $E_q$ be the component of the exceptional divisor for \[\Bl_Q(T\times T) \ra T\times T\] over $q \in Q$. We write $\olo_M(-E_q)$ for the pullback of $\olo(-E_q)$ to $M$. The rank 2 vector bundle $\olo_M(-E_q)\oplus \olo_M$ is Mumford-destabilised by the line bundle $\olo_M(-E_q)$. Fix a K\"ahler class $\Omega_M$ on $M$. By Lemma \ref{bundles} the projective bundle \[\pi: \proj(\olo(-E_q)\oplus \olo_M) \ra M\] admits no cscK metric in the classes $\olo_{\proj}(1) + r\pi^*\Omega_M$ for $r \gg 0$.

Note that by taking a projective bundle we have changed to homotopy type of $M$ so Voisin's Theorem does not immediately imply that $\proj(\olo(-E_q)\oplus \olo_M)$ does not have projective deformations. However this is guaranteed by the following result on deformations of projective bundles \cite{dem_tori}.
\begin{theorem}[Demailly-Eckl-Peternell] Let $X$ be a compact complex manifold with a deformation $X \cong \mathcal{X}_0 \hookrightarrow \mathcal{X} \ra \Delta$. If $X = \proj(E)$ for some vector bundle $E \ra Y$ then $\mathcal{X}_t = \proj(V_t)$ for some deformation $\mathcal{V}_t \ra \mathcal{Y}_t$. 
\end{theorem}
\begin{corollary}\label{non_proj_example}
The K\"ahler manifold $X = \proj(\olo_M(-E_q)\oplus \olo_M)$ is slope unstable (and so has no cscK metrics) with respect to the classes 
\[c_1(\olo_{\proj}(1)) + r\pi^*\Omega_M\]
for $r \gg 0$. Moreover $X$ has no projective deformations.
\end{corollary}
\bibliographystyle{amsplain}

\vskip.3cm
\noindent Universit\`a di Pavia, Via Ferrata 1 27100 Pavia, ITALY\\
and\\
Imperial College, London SW7 2AZ, UK.\\
\textit{E-mail:} jacopo.stoppa\texttt{@unipv.it}
\end{document}